\documentclass[12pt,reqno]{amsart}
\usepackage{amssymb,latexsym}
\usepackage{amsmath}
\usepackage{enumerate}
\usepackage{psfrag}
\usepackage{epsf}
\usepackage{epsfig}

\makeatletter
\@namedef{subjclassname@2010}{
  \textup{2010} Mathematics Subject Classification}
\makeatother

\newtheorem{thm}{Theorem}

\newtheorem{cor}[thm]{Corollary}
\newtheorem{lem}[thm]{Lemma}

\newtheorem{prop}[thm]{Proposition}
\newtheorem{conj}[thm]{Conjecture}

\theoremstyle{definition}

\newcommand{\N}{\mathbb{N}}

\newcommand{\C}{\mathbb{C}}
\newcommand{\Z}{\mathbb{Z}}

\newcommand{\bi}{\begin{itemize}}
\newcommand{\ei}{\end{itemize}}
\newcommand{\bee}{\begin{enumerate}}\newcommand{\ee}{\end{enumerate}}
\newcommand{\bt}{\begin{thm}}\newcommand{\et}{\end{thm}}
\newcommand{\bl}{\begin{lemma}}\newcommand{\el}{\end{lemma}}
\newcommand{\bc}{\begin{cor}}\newcommand{\ec}{\end{cor}}
\newcommand{\bepn}{\begin{prop}}
\newcommand{\enpn}{\end{prop}}
\newcommand{\bepr}{\begin{proof}}
\newcommand{\enpr}{\end{proof}}
\newcommand{\beq}{\begin{equation}}
\newcommand{\eeq}{\end{equation}}
\def\R{{\mathbb R}}
\def\Q{{\mathbb Q}}
\def\Z{{\mathbb Z}}
\def\N{{\mathbb N}}

\def\t{\tau}\def\om{\omega}

\def\Gal{\operatorname{Gal}}
\def\al{\alpha}
\def\la{\lambda}
\def\s{\sigma}

\def\SL{\operatorname{SL}}\def\Gal{\operatorname{Gal}}
\def\PSL{\operatorname{PSL}}
\def\Re{\operatorname{Re}}\def\Im{\operatorname{Im}}

\def\bei{\begin{itemize}}\def\eni{\end{itemize}}
\def\bee{\begin{equation}}\def\ene{\end{equation}}

 \def\geq{\geqslant}

  \def\eps{\varepsilon}

\begin{document}

\baselineskip=17pt

\title[Salem numbers]{  SURVEY ARTICLE: Seventy years of Salem numbers}

\author[C. J. Smyth]{Chris  Smyth}
\address{School of Mathematics and Maxwell Institute for Mathematical Sciences\\University of Edinburgh\\Edinburgh, EH9 3JZ, Scotland}
\email{c.smyth@ed.ac.uk}

\date{}

\begin{abstract}
I survey results about, and recent applications of, Salem numbers.
\end{abstract}



\maketitle

\section{Introduction} In this article I state and prove some basic results about Salem numbers, and then survey
some of the literature about them. My intention is to complement
other general treatises on these numbers, rather than to repeat their
coverage. This applies particularly to the work of Bertin and her
coauthors \cite{Bsurvey,Bbook,Blectnotes} and to the application-rich
Salem number survey of Ghate and Hironaka \cite{GH}. I have, however, quoted some fundamental results from Salem's classical monograph \cite{Sa}.

Recall that a complex number is an {\it algebraic integer} if it is the zero of a polynomial with integer coefficients and leading coefficient $1$. Then its (Galois) {\it conjugates} are the zeros of its {\it minimal polynomial}, which is the lowest degree polynomial of that type that it satisfies. This degree is the {\it degree} of the algebraic integer.

A {\it Salem number} is a real algebraic integer $\t>1$ of degree at
least $4$, conjugate to $\t^{-1}$, all of whose conjugates,
excluding $\t$ and $\t^{-1}$, have modulus $1$. Then $\t+\t^{-1}$ is a
real algebraic integer $>2$, all of whose conjugates $\ne
\t+\t^{-1}$lie in the real interval $(-2,2)$. Such numbers are easy
to find: an example is $\t+\t^{-1}=\frac12(3+\sqrt{5})$, giving
$(\t+\t^{-1}-\frac32)^2=\frac54$, so that $\t^4-3\t^3+3\t^2-3\t+1=0$ and
$\t=2.1537\dots$. We note that this polynomial is a so-called
(self)-reciprocal polynomial: it satisfies the equation $z^{\deg
P}P(z^{-1})=P(z)$. This means that its coefficients form a
palindromic sequence: they read the same backwards as forwards. This
holds for the minimal polynomial of every Salem number. It is simply
a consequence of $\t$ and $\t^{-1}$ having the same minimal
polynomial. Salem numbers are named after  Rapha\" el Salem, who, in 1945, first defined and studied them \cite[Section 6]{S45}.

Salem numbers are usually defined in an apparently more general way,
as in the following lemma. It shows that this
apparent greater generality is illusory.

\begin{lem}[{{Salem \cite[p.26]{Sa}}}] Suppose that $\t>1$ is a real algebraic integer, all of whose conjugates $\ne\t$ lie in the closed unit disc $|z|\le 1$, with at least one on its boundary $|z|=1$. Then $\t$ is a Salem number (as defined above).
\end{lem}
\bepr Taking $\t'$ to be  a conjugate of $\t$ on $|z|=1$, we have that $\overline{\t'}=\t'^{-1}$ is also a conjugate $\t''$ say, so that $\t'^{-1}=\t''$. For any other conjugate $\t_1$ of $\t$ we can apply a Galois automorphism mapping $\t''\mapsto \t_1$ to deduce that $\t_1=\t_2^{-1}$ for some conjugate $\t_2$ of $\t$. Hence the conjugates of $\t$ occur in pairs $\t',\t'^{-1}$. Since $\t$ itself is the only conjugate in $|z|>1$, it follows that $\t^{-1}$ is the only conjugate in $|z|<1$, and so all conjugates of $\t$ apart from $\t$ and $\t^{-1}$ in fact lie on $|z|=1$.
\enpr

It is known that an algebraic integer lying with all its conjugates on the unit circle must be a root of unity (Kronecker \cite{Kr}). So in some sense Salem numbers are the algebraic integers that are  `the nearest things to roots of unity'.  And, like roots of unity, the set of all Salem numbers is closed under taking powers.

\begin{lem}[{{Salem \cite[p.169]{S45}}}]\label{L-2} If $\t$ is a Salem number of degree $d$, then so is $\t^n$ for all $n\in \N$.
\end{lem}
\bepr If $\t$ is conjugate to $\t'$ then $\t^n$ is conjugate to $\t'^n$. So $\t^n$ will be a Salem number of degree $d$ unless some of its conjugates coincide: say $\t_1^n=\t_2^n$ with $\t_1\ne\t_2$. But then, by applying a Galois automorphism mapping $\t_1\mapsto \t$, we would have $\t^n=\t_3^n$ say, where $\t_3\ne \t$ is a conjugate of $\t$, giving $|\t^n|>1$ while $|\t_3^n|\le 1$, a contradiction.
\enpr
Which number fields contain Salem numbers? Of course one can simply choose a list of Salem numbers $\t,\t',\t'',\dots$ say, and then the number field $\Q(\t,\t',\t'',\dots)$ certainly contains $\t,\t',\t'',\dots$. However, if one is interested only in finding all Salem numbers in a field $\Q(\t)$ for $\t$ a given Salem number, we can be much more specific.

\begin{prop}[{{Salem \cite[p.169]{S45}}}]\label{P-3} \bei\item[(i)] A number field $K$ is of the form $\Q(\t)$ for some Salem number $\t$ if and only if $K$ has a totally real subfield $\Q(\al)$ of index $2$, and $K=\Q(\t)$ with $\t+\t^{-1}=\al$, where $\al>2$ is an irrational algebraic integer, all of whose conjugates $\ne\al$ lie in $(-2,2)$.

 \item[(ii)] Suppose that $\t$ and $\t'$ are Salem numbers with $\t'\in\Q(\t)$. Then $\Q(\t')=\Q(\t)$ and, if $\t'>\t$, then $\t'/\t$ is also a Salem number in $\Q(\t)$.

\item[(iii)] If $K=\Q(\t)$ for some Salem number $\t$, then there is a Salem number $\t_1\in K$ such that the set of Salem numbers  in $K$ consists of the powers of $\t_1$.
\eni
\end{prop}
\bepr \bei\item[(i)] Suppose that $K=\Q(\t)$ for some Salem number $\t$. Then, defining $\al=\t+\t^{-1}$, every conjugate $\al'$ of $\al$ is of the form $\al_1=\t_1+\t_1^{-1}$, where $\t_1$, a conjugaste of $\t$, either is one of $\t^{\pm 1}$ or has modulus $1$. Hence $\al_1$ is real, and so $\al$ is totally real with $\al>2$ and all other conjugates of $\al$ in $(-2,2)$. Also, because $\Q(\al)$ is a proper subfield of $\Q(\t)$ and $\t^2-\al\t+1=0$ we have $[\Q(\t):\Q(\al)]=2$.

Conversely, suppose that $K$ has a totally real subfield $\Q(\al)$ of index $2$, and $K=\Q(\t)$, where $\al>2$ is an irrational algebraic integer, all of whose conjugates $\ne\al$ lie in $(-2,2)$, and  
$\t+\t^{-1}=\al$. Then  for every conjugate $\t_1$ of $\t$ we have $\t_1+\t_1^{-1}=\al_1$ for some
conjugate $\al_1$ of $\al$, and for every conjugate $\al_1$ of $\al$ we have $\t_1+\t_1^{-1}=\al_1$ for some
conjugate $\t_1$ of $\t$. If $\al_1=\al$ then $\t_1=\t^{\pm 1}$. Otherwise, $\al_1\in(-2,2)$ and so $|\t_1|=1$. Furthermore, as $\al$ is irrational it does indeed have a conjugate $\al_1\in(-2,2)$ for which the corresponding $\t_1$ has modulus $1$. Hence $\t$ is a Salem number.

\item[(ii)] By  (i) we can write $\t'$ as $\t'=p(\al)+\t q(\al)$, where $p(z),q(z)\in\Q[z]$ and $\al=\t+\t^{-1}$. If $\t'$ had lower degree than $\t$, say $[\Q(\t):\Q(\t')]=k>1$ then, for the $d=[\Q(\t):\Q]$ conjugates $\t_i$ of $\t$, $k$ of the values $p(\t_i+\t_i^{-1})+\t_i q(\t_i+\t_i^{-1})$ would equal $\t'$ and another $k$ of them would equal ${\t'}^{-1}$. In particular, $p(\t_i+\t_i^{-1})+\t_i q(\t_i+\t_i^{-1})$ would be real for some nonreal $\t_i$, giving $q(\t_i+\t_i^{-1})=0~$ and hence, on applying a suitable automorphism, that $q(\al)=0$. Thus we would have $\t'=p(\al)$, a contradiction, since $p(\al)$ is totally real. So $k=1$ and $\Q(\t')=\Q(\t)$.

 Since $\t'\in \Q(\t)$, it is a polynomial in $\t$. Therefore any Galois automorphism taking $\t\mapsto\t^{-1}$ will map $\t'$ to a real conjugate of $\t'$, namely $\t'^{\pm 1}$. But it cannot map $\t'$ to itself for then, as $\t$ is also a polynomial in $\t'$, $\t$ would be mapped to itself, a contradiction. So $\t'$ is mapped to $\t'^{-1}$ by this automorphism. Hence $\t'\t^{-1}$ is conjugate to its reciprocal.  So the conjugates of $\t'\t^{-1}$ occur in pairs $x, x^{-1}$. Again,  because $\t'$ is a polynomial in $\t$, any automorphism fixing $\t$ will also fix $\t'$, and so fix $\t'\t^{-1}$. Likewise, any automorphism fixing $\t'$ will also fix $\t$.

Next consider any conjugate of $\t'\t^{-1}$ in $|z|>1$. It will be of the form
$\t_1'\t_1^{-1}$, where $\t_1'$ is a conjugate of $\t'$ and $\t_1$ is a conjugate of $\t$. For this to lie in  $|z|>1$, we must either have $|\t_1'|>1$ or $|\t_1|<1$, i.e., $\t_1'=\t'$ or $\t_1=\t^{-1}$. But in the first case, as we have seen, $\t_1=\t$, so that $\t_1'\t_1^{-1}=\t'\t^{-1}$, while in the second case $\t_1'=\t'^{-1}$, giving $\t_1'\t_1^{-1}=\t\t'^{-1}<1$. Hence
 $\t'\t^{-1}$ itself is the only conjugate of $\t'\t^{-1}$ in $|z|>1$. It follows that all conjugates of $\t'\t^{-1}$ apart from $(\t'\t^{-1})^{\pm 1}$ must lie on $|z|=1$, making $\t'\t^{-1}$ a Salem number.

\item[(iii)] Consider the set of all Salem numbers in $K=\Q(\t)$. Now the number of Salem numbers $<\t$ in $K$ is clearly finite, as there are only finitely many possibilities for the minimal polynomials of such numbers. Hence there is a smallest such number,  $\t_1$ say. For any  Salem number, $\t'$ say, in $K$ we can choose a positive integer $r$ such that $\t_1^r\le \t'<\t_1^{r+1}$. But if $\t_1^r< \t'$ then, by (ii) and Lemma \ref{L-2},  $\t'\t_1^{-r}$ would be another Salem number in $K$ which, moreover, would be less than $\t_1$, a contradiction.
Hence $\t'=\t_1^r$.
\eni

\enpr

The statements (ii) and (iii) above in fact represent a slight strengthening of Salem's results, as they  do not assume that all Salem numbers in $K$ have degree $[K:\Q]$.

We now show that the powers of Salem numbers have an unusual property.

\begin{prop} [{{Essentially Salem \cite{S45}, reproduced in \cite{Sa}}}]\label{P-Teps} For every Salem number $\t$ and every $\eps>0$ there is a  real number $\la>0$ such that  the distance $\|\la\t^n\|$ of $\la\t^n$ to the nearest integer is less than $\eps$ for all $n\in\N$.
\end{prop}

The result's first explicit appearance seems to be in Boyd \cite{B-bad}, who credits Salem. A proof comes from an easy modification the proof in \cite[pp 164-166]{S45}.
\bepr We consider the standard embedding of the algebraic integers  $\Z(\t)$ as a lattice in $\R^d$ defined  for $k=0,1,\dots,d-1$ by the map
\[
\t^k\mapsto(\t^k,\t^{-k},\Re\t_2^k,\Im\t_2^k,\Re\t_3^k,\Im\t_3^k,\dots,\Re\t_{d/2}^k,\Im\t_{d/2}^k),
\]
 where $\t^{\pm 1},\t_j^{\pm 1}\,(j=2,\dots,d/2)$ are the conjugates of $\t$. As this is a lattice of full dimension $d$, we know that for every $\eps'>0$ there are lattice points in the `slice' $\{(x_1,\dots,x_n)\in\R^n\, :\, |x_i|<\eps'\,(i=2,\dots,d)\}$. Such a lattice point corresponds to an element $\la(\t)$ of $\Z(\t)$ with conjugates $\la_i$ satisfying $|\la_i|<\sqrt{2}\eps'\,(i=2,\dots,d)$.

Next, consider the sums
\[
\s_n=\la(\t)\t^n+\la(\t^{-1})\t^{-n}+\la(\t_2)\t_2^n+\la(\t_2^{-1})\t_2^{-n}+\dots \la(\t_{d/2})\t_{d/2}^n+\la(\t_{d/2}^{-1})\t_{d/2}^{-n},
\]
where $\la(x)\in\Z[x]$. Since $\s_n$ is a symmetric function of the conjugates of $\t$, it is rational. As it is an algebraic integer, it is in fact a rational integer. Since all terms $\la(\t^{-1})\t^{-n}$, $\la(\t_2)\t_2^n$, $\la(\t_2^{-1})\t_2^{-n}$, \dots,  $\la(\t_{d/2})\t_{d/2}^n$, $\la(\t_{d/2}^{-1})\t_{d/2}^{-n}$ are $<\sqrt{2}\eps'$ in modulus, we see that
\[
|\s_n-\la(\t)\t^n|<(d-1)\sqrt{2}\eps'.
\]
Hence, choosing $\eps'=\eps/((d-1)\sqrt{2})$, we have $\|\la\t^n\|\le|\s_n-\la(\t)\t^n|<\eps$.
\enpr

In fact, this property essentially characterises Salem (and Pisot) numbers among all real numbers. Recall that a {\it Pisot number} is an algebraic integer greater than $1$ all of whose conjugates, excluding itself,  all lie in the open unit disc $|z|<1$. Pisot \cite{Pi} proved that if  $\la$ and $\t$ are real numbers such that
\bee\label{EPi}
\|\la\t^n\|\le \frac{1}{2e\t(\t+1)(1+\log \la)}\qquad (= B \text{ say}),
\ene
for all integers $n\ge 0$ then $\t$ is either a Salem number or a Pisot number and $\la\in\Q(\t)$.
The denominator in this result was later improved by Cantor \cite{Ca}  to $2e\t(\t+1)(2+\sqrt{\log \la})$, and then by Decomps-Guilloux and Grandet-Hugot \cite{DG}  to
$e(\t+1)^2(2+\sqrt{\log \la})$.
However, Vijayaraghavan \cite{V} proved that for each $\eps>0$ there are uncountably many real numbers $\al>1$ such that $\|\al^n\|<\eps$ for all $n\ge 0$. To be compatible with \eqref{EPi}, it is clear that such $\al$ that are not Pisot or Salem numbers must be large (depending on $\eps$). Specifically, if $\al>(2e\eps)^{-1/2}$ then there is no contradiction to \eqref{EPi}.  Furthermore, Boyd \cite{B-bad} proved that  if the bound $B$ in \eqref{EPi} is replaced by $10B$ then $\t$ can be transcendental.

For more results concerning the distribution of the fractional parts of $\la\t^n$ for $\t$ a Salem number, see Dubickas \cite{D2}, Za\"\i mi \cite{Z1,Z2,Z3}, and Bugeaud's monograph \cite[Section 2.4]{Bu}.

\section{A smallest Salem number?} Define the polynomial $L(z)$ by
\[
L(z)=z^{10}+z^9-z^7-z^6-z^5-z^4-z^3+z+1.
\]
This is the minimal polynomial of the Salem number $\t_{10}=1.176\dots$, discovered by
 D. H. Lehmer \cite{L1} in 1933 (i.e., before Salem numbers had been defined!). Curiously, the polynomial $L(-z)$ had appeared a year
 earlier in Reidemeister's book \cite{Re} as the Alexander polynomial of the $(-2,3,7)$
 pretzel knot.  Lehmer's paper seems to be the first where what is now called the {\it Mahler measure} of a polynomial
 appears: the Mahler measure $M(P)$ of a monic one-variable polynomial $P$ is the product
 $\prod_i\max(1,|\al_i|)$ over the roots $\al_i$ of the polynomial. For a survey of Mahler measure see \cite{Ssurv}; for  its multivariable generalisation, see   Boyd \cite{Boyd-spec} and
 Bertin and Lal\'\i n \cite{Bertin-Lalin}.

  Lehmer also asked whether the  Mahler measure of any nonzero noncyclotomic
 irreducible polynomial with integer coefficients is bounded below by some constant $c>1$.
  This is now commonly referred to as `Lehmer's problem' or (inaccurately) as `Lehmer's conjecture' --- see \cite{Ssurv}. If indeed there were such a bound,
  then certainly Salem numbers would be bounded away from $1$, but this would not immediately imply that there is a smallest Salem number (but see the end of Section \ref{SS-3.1}). However, the `strong version' of `Lehmer's conjecture'
  states that one can take $c=\t_{10}$, implying that there is indeed a smallest Salem number, namely $\t_{10}$. A consequence of this strong version is the following.

\begin{conj}\label{CS} Suppose that $d\in\N$ and $\al_1,\al_2,\dots,\al_d$ are real numbers with
$\al_1\in(2,\t_{10}+\t_{10}^{-1})$ and $\al_2,\al_3,\dots,\al_d\in(-2,2)$.
Then $\prod_{i=1}^d(x+\al_i)\not\in\Z[x]$.
\end{conj}
(Note that $\t_{10}+\t_{10}^{-1}=2.026\dots$.)  For if there were $\al_1,\al_2,\dots,\al_d$ in the  intervals stated with $\prod_{i=1}^d(x+\al_i)\in\Z[x]$, then the algebraic integer $\t>1$ defined by $\t+\t^{-1}=\al_1$ would be a Salem number less than $\t_{10}$.

\section{Construction of Salem numbers}
\subsection{Salem's method}\label{SS-3.1}
Salem \cite[Theorem IV, p.30]{Sa} found a simple way to construct infinite sequences of Salem numbers from Pisot numbers.  Now if $P(z)$ is the minimal polynomial of a Pisot number, then, except possibly for some small values of $n$, the polynomials $S_{n,P,\pm 1}(z)=z^nP(z)\pm z^{\deg P}P(z^{-1})$ factor as the minimal polynomial of a Salem number, possibly multiplied by some cyclotomic polynomials. In particular, for $P(z)=z^3-z-1$, the minimal polynomial of the smallest Pisot number, $S_{8,P,-1}=(z-1)L(z)$. Salem's construction shows that every Pisot number is the limit on both sides of a sequence of Salem numbers. (The construction has to be modified slightly when $P$ is reciprocal; this occurs only for certain $P$ of degree $2$.)

Boyd \cite{B1} proved that all Salem numbers could be produced by Salem's construction, in fact with $n=1$. It turns out that many different Pisot numbers can be used to produce the same Salem number. These Pisot numbers can be much larger than the Salem number they produce.
In particular, on taking $P(z)=z^3-z-1$ and $\eps=-1$, the minimal polynomial of the smallest Pisot number $\theta_0=1.3247\dots$, Salem's method shows that there are infinitely many Salem numbers less than $\theta_0$. This fact motivates the next definition, due to Boyd.

    Salem numbers less than $1.3$ are called {\it small}. A list of $39$ such numbers was compiled by Boyd \cite{B1}, with later additions of four each by Boyd \cite{B2} and Mossinghoff \cite{Mo1}, making $47$ in all -- see \cite{Mo2}. (The starred entries in this list are the four Salem numbers found by Mossinghoff. They include one of degree $46$.) Further, it was determined by Flammang,  Grandcolas and Rhin \cite{FGR} that the  table was complete up to degree $40$.  This was extended up to degree $44$ by Mossinghoff, Rhin and Wu \cite{MRW} as part of a larger project to find small Mahler measures. Their result shows that if Conjecture \ref{CS} is false then any counterexample to that conjecture must have degree $d\ge 23$.

 In \cite{B3} Boyd showed how to find, for a given $n\ge 2$, $\eps=\pm 1$ and real interval $[a,b]$, all Salem numbers in that interval that are roots of $S_{n,P,\eps}(z)=0$ for some Pisot number having minimal polynomial $P(z)$. In particular, of the four new small Salem numbers that he found,  two were discovered by this method. The other two he found in \cite{B3} are not of this form: they are roots only of some $S_{1,P,\eps}(z)=0$.

 Boyd and Bertin \cite{BB0,BB} investigated the properties of the polynomials $S_{1,P,\pm 1}(z)$ in detail. For a related, but interestingly different, way of constructing Salem numbers, see Boyd and Parry \cite{BP}.

Let $T$ denote the set of all Salem numbers (Salem's notation). (It couldn't be called $S$, because that is used for the set of all Pisot numbers. The notation $S$ here is in honour of Salem, however: Salem \cite{Sclosed} had proved the magnificent result that the Pisot numbers form a closed subset of the real line.) Salem's construction shows that the derived set (set of limit points) of $T$ contains $S$. Salem \cite[p.31]{Sa} wrote `We do not know whether numbers of $T$ have limit points other than $S$'. Boyd \cite[p. 327]{B1} conjectured that  there were no other such limit points, i.e., that the derived set of $S\cup T$ is $S$.
 (Not long before, he had  conjectured \cite{B15} that $S\cup T$ is closed -- a conjecture that left open the possibility that some numbers in $T$ could be limit points of $T$.)

Salem's construction shows that every Pisot number is a limit from below of Salem numbers. So if the derived set of $S\cup T$ is indeed $S$, then any limit point of Salem numbers from above is also a limit point of Salem numbers from below. Hence Boyd's conjecture implies that there must be a smallest Salem number.

\subsection{Salem numbers and matrices}
One strategy that has been used to try to solve Lehmer's Problem is to attach some combinatorial object (knot, graph, matrix,\dots) to an algebraic number (for example, to a Salem number). But it is not clear whether
the object could throw light on the (e.g.) Salem number, or, on the contrary, that the Salem number could throw light on the object.

Typically, however, such attachment constructions seem to  work only for a restricted class of algebraic numbers, and not in full generality. For example, McKee, Rowlinson, and Smyth   considered star-like trees as the objects for attachment. this was extended by McKee and Smyth to  more general graphs \cite{McKSm3} and then to integer symmetric matrices \cite{McKSm-ISM,McKSm-IMRN}.  (These can be considered as generalisations of graphs: one can identify a graph with its adjacency matrix -- an integer symmetric matrix having all entries $0$ or $1$, with only zeros on the diagonal.) The main tool for their work was the following classical result, which deserves to be better known.

\begin{thm} [{{Cauchy's Interlacing Theorem -- for a proof see for instance \cite{Fisk}}}] Let $M$ be a real $n\times n$ symmetric matrix, and $M'$ be the matrix obtained from $M$ by removing the $i$th row and column. Then the eigenvalues $\la_1,\dots,\la_n$ of $M$ and the eigenvalues $\mu_1,\dots,\mu_{n-1}$ of $M'$ interlace, i.e.,
\[
\la_1\le\mu_1\le\la_2\le\mu_2\le\dots\le\mu_{n-1}\le\la_n.
\]
\end{thm}

We say that an $n\times n$ integer symmetric matrix $M$ is {\it cyclotomic} if all its eigenvalues lie in the interval $[-2,2]$. It is so-called because then its associated reciprocal polynomial
\[
 P_M(z)=z^n\det\left((z+z^{-1})I-M\right)
\]
 has all its roots on $|z|=1$ and so (Kronecker again) is a product of cyclotomic
 polynomials. Here $I$ is the $n\times n$ identity matrix.

\begin{figure}[h]
\begin{center}
\leavevmode
\psfragscanon
\psfrag{A}[l]{$\tilde A_n$}
\psfrag{D}[l]{$\tilde D_n$}
\psfrag{6}[l]{$\tilde E_6$}
\psfrag{7}{$\tilde E_7$}
\psfrag{8}{$\tilde E_8$}
\includegraphics[scale=0.4]{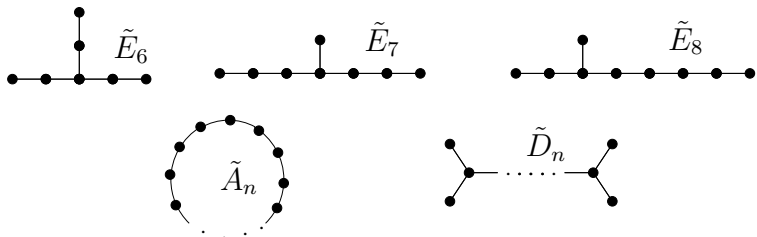}
\end{center}
\caption{The Coxeter graphs $\tilde E_6, \tilde E_7,
 \tilde E_8, \tilde A_n (n\geq 2)$ and $ \tilde
D_n (n\geq 4)$. (The number of vertices is $1$ more than the index.)}\label{F-except}
\end{figure}

The cyclotomic graphs are very familiar.
\bt [{{ J.H. Smith  \cite{Smi}}}] The connected cyclotomic graphs consist of the (not necessarily proper) induced subgraphs of the Coxeter graphs $\tilde A_n(n\ge 2)$, $\tilde D_n(n\ge 4)$, $\tilde E_6$, $\tilde E_7$, $\tilde E_8$, as in Figure \ref{F-except}.
\et

(These graphs also occur in the theory of  Lie algebras, reflection groups, Lie groups, Tits geometries, surface singularities, subgroups of $\operatorname{SU}_2(\C)$ (McKay correspondence),\dots).

McKee and Smyth \cite{McKSm-ISM} describe all the cyclotomic matrices, of which the cyclotomic graphs form a small subset. They prove in \cite{McKSm-IMRN} that the strong version of Lehmer's conjecture is true for the set of polynomials $P_M$: namely, if $M$ is not a cyclotomic matrix, then $P_M$ has Mahler measure at least $\t_{10}=1.176\dots$, the smallest known Salem number. In fact they show that the smallest three known Salem numbers are all Mahler measures of $P_M$ for some integer symmetric matrix $M$, while the fourth-smallest known Salem number is not.

Most constructions of families of Salem numbers produce a monic reciprocal integer polynomials having roots $\t>1$, $1/\t<1$ and all other roots on the unit circle. The final stage of the construction requires that each polynomial be expressed as a product of irreducible factors. Then, by Kronecker's Theorem \cite{Kr}, the irreducible factor having $1/\t$ as a root must also have $\t$ as a root, and any other factors must be cyclotomic polynomials. For $\t$ to be a Salem number, the first-mentioned factor must also have a root of modulus $1$. To determine its degree, it is usually necessary to determine the degrees of the cyclotomic factors. One method, used by Smyth \cite{Sm2} and  McKee and Yatsyna \cite{McKY} (see also  Beukers and Smyth \cite{BeSm}), is to make use of the fact \cite[Lemma 2.1]{Sm2} that every root of unity $\om$ is conjugate to one of $\om^2$, $-\om^2$ or $-\om$. A second method is to use results of Mann \cite{Mann} on sums of roots of unity. This method was used by  Gross, Hironaka and McMullen \cite{GHM} for polynomials coming from modified Coxeter diagrams, and by Brunotte and Thuswaldner \cite{BT} for polynomials coming from star-like trees.

For other construction methods for Salem numbers see  Lakatos \cite{L1,L2,L3} and also \cite{McKSm1,McKSm2,McKSm3,McKSm4,MRS,Sm1,Sm2}. In particular, in \cite{L1,L3} Lakatos shows that Salem numbers arise as the spectral radius of Coxeter transformations of certain oriented graphs containing no oriented cycles.

\subsection{Traces of Salem numbers}
McMullen \cite[p.230]{mcm1} asked whether there are any Salem numbers of trace less than $-1$. McKee and Smyth \cite{McKSm1,McKSm2} found examples of Salem numbers of trace $-2$, and indeed showed that there are Salem numbers of every trace. 
It is known \cite{Sm2} that the smallest degree of a Salem number of trace $-1$ is $8$, and there are Salem numbers of trace $-1$ for every even degree at least $8$.  Recently McKee and Yatsyna \cite{McKY} have shown that there are Salem numbers of trace $-2$ for every even degree at least  $38$.

 If $\t$ has degree $d\ge 10$ then its  trace is known to be at least at least $\lfloor 1-d/9\rfloor$ -- see \cite{McKSm2}. Conversely, Salem numbers having  trace $-T\le -2$ must have degree $d\ge 2\lfloor 1+\frac92 T\rfloor$ (corrected from \cite[p. 35]{McKSm2}). For $-T=-2$ this bound is attained, it having been shown in \cite{McKSm1} that there are two Salem numbers of trace $-2$ and degree $20$.
For a summary of the status of even degrees $20<d<38$ for which there are known to be or not to be  Salem numbers of degree $d$ and trace $-2$, see \cite[Section 3]{McKY}.

 For $-T=-3$ the above bound gives $d\ge 28$,  later improved  by Flammang \cite{Fla}  to $d\ge 30$. In the other direction  El Omani, Rhin and Sac-\' Ep\' ee \cite{EORSE} have recently found Salem number examples of trace $-3$ and degree $34$. Their computations strongly suggest that $34$ is in fact the smallest degree for Salem numbers of trace $-3$.
Recently Liang and Wu \cite{LW} have shown that a Salem number of trace $-4$ must have degree at least $40$, and a Salem number of trace $-5$ must have degree at least $50$.

All of these lower bounds for the degree of a Salem number of given negative trace  make use of the fact that, for a Salem number $\t$, the number $\t+1/\t+2$ is a totally positive algebraic integer. Thus known lower bounds of the type $\text{trace}(\alpha)>\la d'$ for totally positive algebraic integers
 $\al$ of degree $d'\ge 5$ can be used. Specifically, one readily deduces that, from such an inequality, a Salem number of degree $d\ge 10$ has trace at least $\lfloor 1-(1-\la/2)d\rfloor$, and a Salem number of trace $-T\le -2$ has degree at least $2\lfloor\frac{T}{2-\la}\rfloor+2$. (The results quoted above are for $\la=16/9$.)

 For an interesting survey of the trace problem for totally positive algebraic integers see  Aguirre and Peral \cite{AP}.

\subsection{Distribution modulo 1 of the  powers of a Salem number}\label{S-dense}
Let $\t>1$ be a Salem number of degree $d$. Salem \cite[Theorem V, p.33]{Sa} proved that although the powers $\t^n\,\text{mod }1$ of $\t$ are everywhere dense on $(0,1)$, they are not uniformly distributed on this interval. In fact, the asymptotic distribution function of $\t^n\,\text{mod }1$ is that of a sum $2\sum_{j=1}^{d/2-1}\cos{\boldsymbol\theta}_j$, where the ${\boldsymbol\theta}_j$ are independent random variables uniformly distributed on $[0,\pi]$. In the simplest case, $d=4$, a routine calculation shows that for $x\in[0,1]$
\begin{align*}
\text{Prob}(&\t^n\,\text{mod }1\le x)=\\ &\tfrac12+\frac{\arcsin(\frac{x+1}{2})+\arcsin(\frac{x}{2})+\arcsin(\frac{x-1}{2})+\arcsin(\frac{x-2}{2})}{\pi}.
\end{align*}
In particular,  $\text{Prob}(1/4\le \t^n\,\text{mod }1\le 3/4)=0.4134038362$, significantly less than $1/2$.

More generally,  Akiyama and Tanigawa \cite{AT} gave a quantitative description of how far the sequence $\t^n\,\text{mod }1$ is from being uniformly distributed. They show,
for $\tau$ a Salem number of degree $2d'\ge 8$ and $A_N(\{\tau^n\},I)$ being the number of $n\le N$ for which the fractional part $\{\t^N\}$ lies in a subinterval $I$ of $[0,1]$,
that $\lim_{N\to\infty} \frac{1}{N}A_N(\{\tau^n\},I)$ exists and satisfies
\[
\left|\lim_{N\to\infty} \tfrac{1}{N}A_N(\{\tau^n\},I)-|I|\right|\le 2\zeta\left(\tfrac12(d'-1)\right)(2\pi)^{1-d'}|I|.
\]
Here  $|I|$ is the length of $I$.
Note that this difference tends to $0$ as $d'\to\infty$, so that $\{\tau^n\}$ is more nearly uniformly distributed mod $1$ for $\t$ of large degree.
See also \cite{DMR}.

\subsection{Sumsets of Salem numbers} Dubickas \cite{D1} shows that a sum of $m\ge 2$ Salem numbers cannot be a Salem number, but that for every $m\ge 2$ there are $m$ Salem numbers whose sum is a Pisot number and also $m$ Pisot numbers whose sum is a Salem number.

\subsection{Galois group of Salem number fields} Lalande \cite{Lal} and Christopoulos and McKee \cite{CMcK}
studied the Galois group of a number field defined by a Salem
number.  Let $\t$ be a Salem number of degree $2n$, $K=\Q(\t)$ and
$L$ be its Galois closure. Then it is known that $G:=\Gal(L/\Q)\le
C_2^n\rtimes S_n$. Conversely, if $K$ is a real number field of
degree $2n>2$ with exactly $2$ real embeddings, and, for its Galois
closure $L$, that $G\le C_2^n\rtimes S_n$, then Lalande proved that
$K$ is generated by a Salem  number.

Now, for a Salem number $\t$, let $K'=\Q(\t+\t^{-1})$, $L'$ be its
Galois closure and $N\subset G$ be the fixing group of $L'$. Then
Christopoulos and McKee showed that $G$ is isomorphic to $N\rtimes
\Gal(L'/\Q)$, where $N$ is isomorphic to either $C_2^n$ or
$C_2^{n-1}$. The latter case is possible only when $n$ is odd.

Amoroso \cite{A} found a lower bound, conditional on the Generalised
Riemann Hypothesis, for the exponent of the class group of such
number fields $L$.

\subsection{The range of polynomials $\Z[\t]$}
 P. Borwein and Hare \cite{BH} studied the `spectrum' of values $a_0+a_1\t+\dots+a_n\t^n$ when the $a_i\in\{-1,1\}$, $n\in\N$ and $\t$ is a Salem number. They showed that if  $\t$ was a Salem number defined by being the zero of a polynomial of the form $z^m-z^{m-1}-z^{m-2}-\cdots-z^2-z+1$ for some $m\ge 4$
, then this spectrum is discrete. They also asked \cite[Section 7]{BH}
\bei
\item Are these the only Salem numbers with this spectrum discrete?
\item  Are the only $\t$ where this spectrum is discrete and $M(\t)<2$ necessarily Salem numbers or Pisot numbers?
\eni

Hare and Mossinghoff \cite{HM} show, given a Salem number $\t<\frac12(1+\sqrt{5})$ of degree at most $20$, that some sum  of distinct powers of $-\t$ is zero, so that $-\t$ satisfies some Newman polynomial.

 Feng \cite{Fe2} remarked that it follows from Garsia \cite[Lemma 1.51]{Ga} that, given a Salem number $\t$ and $m\in\N$ there exists $c>0$ and $k\in\N$ (depending on $\t$ and $m$) such that for each $n\in\N$
there are no nonzero numbers $\xi=\sum_{i=0}^{n-1} a_i\t^i$ with $a_i\in\Z$, $|a_i|\le m$ and $|\xi|<\frac{c}{n^k}$. He asks whether, conversely, if $\t$ is any non-Pisot number in $(1,m+1)$ with this property, then must $\t$ necessarily be a Salem number?

\subsection{Other Salem number studies} Salem \cite{S45},  \cite[p. 35]{Sa} proved that every Salem number is the quotient of two Pisot numbers.

For connections between small Salem numbers and exceptional units, see Silverman \cite{Sil}.

Dubickas and Smyth \cite{DS} showed that any line in $\mathbb C$ containing two nonreal conjugates of a Salem number cannot contain the Salem number itself.

Akiyama and Kwon \cite{AK} constructed Salem numbers satisfying polynomials whose coefficients are nearly constant.

For generalisations of Salem numbers, see Bertin \cite{BK1,BK2}, Cantor \cite{Ca2}, Kerada \cite{Ke}, Meyer \cite{Me}, Samet \cite{Sam}, Schreiber \cite{Schr} and Smyth \cite{Sm1}. Note the correction made to \cite{Sam} in \cite{Sm1}. See also Section \ref{s1} below for $2$-Salem numbers.

\section{Salem numbers outside Number Theory}

The survey of Ghate and Hironaka \cite{GH} contains many applications of Salem numbers, for the period up to 1999. Only a few of the applications they describe are briefly recalled here, in subsections \ref{s1}, \ref{s2} and \ref{s3}. Otherwise, I concentrate on developments since their paper appeared.

For some of these applications, the restriction that Salem numbers should have degree at least $4$ can be dropped: the results also hold for reciprocal Pisot numbers, whose minimal polynomials are $x^2-ax+1$ for $a\in\N$, $n\ge 3$. Some authors include these numbers in the definition of Salem numbers. Accordingly, I will allow  these numbers to be (honorary!) Salem numbers in this section.

 (Note, however, that for $\t$ such a `quadratic Salem number', the fractional parts of the sequence $\{\t^n\}_{n\in\N}$ tend to $1$ as $n\to\infty$, whereas for true Salem numbers this sequence is dense in $(0,1)$, as stated in Section \ref{S-dense}. Furthermore, neither  Proposition \ref{P-3}(ii) nor (iii) would hold, since the proof of Lemma \ref{P-3}(ii) depends on a Salem number having a nonreal conjugate. Indeed, for the example Salem number $\t=2.1537\dots$ in the introduction, the field $\Q(\t)$ contains the Salem numbers $\t^n\,(n=1,2,\dots)$, with $\t$ the smallest Salem number in this field. However, $\Q(\t)$ also contains the `quadratic Salem number' $\t+\t^{-1}=\frac12(3+\sqrt{5})$, which is not a power of $\t$.)

\subsection{Growth of groups}\label{s1} For a group $G$ with finite generating set $S=S^{-1}$, we define its {\it growth series} $F_{G,S}(x)=\sum_{n=0}^\infty a_nx^n$, where $a_n$ is the number of elements of $G$ that can be represented as the product of $n$ elements of $S$, but not by fewer. For certain such groups, $F_{G,S}(x)$ is known to be a rational function. Then expanding $F_{G,S}(x)$ out in partial fractions leads to a closed formula for the $a_n$. See \cite[Section 4]{GH} for a detailed description, including references. See also  \cite{BCT}.

In particular, let $G$ be a Coxeter group generated by reflections in $d\ge 3$ geodesics in the upper half plane, forming a polygon with angles $\pi/p_i\,(i=1,2,\dots,d)$, where $\sum_i\pi/p_i<\pi$. Taking $S$ to be the set of these reflections, it is known (Cannon and Wagreich, Floyd and Plotnick, Parry) that then the denominator of $F_{G,S}(x)$ -- call it $\Delta_{p_1,p_2,\dots,p_d}(x)$ -- is the minimal polynomial of a Salem number, $\tau$ say, possibly multiplied by some cyclotomic polynomials. Then the $a_n$ grow exponentially with growth rate $\lim_{n\to\infty}a_{n+1}/a_n=\tau$. Hironaka \cite{Hi} proved that among all such $\Delta_{p_1,p_2,\dots,p_d}(x)$, the lowest growth rate was achieved for
$\Delta_{p_1,p_2,p_3}(x)$, which is Lehmer's polynomial $L(x)$, with growth rate $\t_{10}=1.176\dots$.

To generalise a bit, define a real $2$-Salem number to be an algebraic integer $\al>1$  which has exactly one conjugate $\al'\ne \al$ outside the closed unit disc, and at least one conjugate on the unit circle. Then all conjugates of $\al$ apart from $\al^{\pm 1}$ and ${\al'}^{\pm 1}$ have modulus $1$.
Zerht and Zerht-Liebensd\" orfer \cite{ZZ}  give examples of infinitely many cocompact Coxeter groups (``Coxeter Garlands'') in $\mathbb H^4$ with the property that their growth function has denominator
\begin{align*}
D_n(z)&=p_n(z)+nq_n(z)\\
&=z^{16}\! -\! 2z^{15}\! +\! z^{14}\! -\! z^{13}\! +\! z^{12}\! -\! z^{10}\! +\! 2z^9\! -\! 2z^8\! +\! 2z^7\! -\! z^6\! +\! z^4\! -\! z^3\! +\! z^2\! -\! 2z\! +\! 1\\
&\qquad+nz(-2z^{14}+z^{12}+z^{10}+z^9+2z^7+z^5+z^4+z^2-2),
\end{align*}
which, if irreducible, would be the minimal polynomial of a $2$-Salem number.

 Umemoto \cite{Um} showed that $D_1(t)$ is irreducible\footnote{In fact, one can show that $D_n(z)$ is irreducible for all $n\ge 1$. A sketch is as follows: comparison with the table \cite{Mo2} shows that neither root of $D_n(z)$ in $|z|>1$ can be a Salem number. Then putting $z=e^{it}$, the fact that $e^{-8it}p_n(e^{it})/q_n(e^{it})$ is real and $>0$ for small $t>0$ shows that $D_n(z)$ has no cyclotomic factors. (Selberg \cite[p. 705]{Sel} remarks that he has always found a sketch of a proof much more informative than a complete proof.)},
 and also produced infinitely many cocompact Coxeter groups whose growth rate is a $2$-Salem number of degree $18$. The growth rate in these examples is the larger of the two $2$-Salem conjugates that are outside the unit circle. This is compatible with a conjecture of Kellerhals and Perren \cite{KP} that the growth rate of a Coxeter group acting on hyperbolic $n$-space should be a Perron number (an algebraic integer $\al$ whose  conjugates different from $\al$ are all of  modulus less than $|\al|$.) This has been verified for $n=3$ for so-called generalised simplex groups by Komori and  Umemoto \cite{KU}.

For some  other recent papers on non-Salem growth rates see \cite{Kel}, \cite{KK}, \cite{Ko}.

\subsection{Alexander Polynomials}\label{s2} A result of Seifert tells us that a polynomial $P\in\Z[x]$ is the Alexander polynomial of some knot iff it is monic and reciprocal, and $P(1)=\pm 1$. In particular, Hironaka \cite{Hi} showed that
$\Delta_{p_1,p_2,\dots,p_d}(-x)$ is the Alexander polynomial of the $(p_1,p_2,\dots,p_d,-1)$ pretzel knot. Hence, from the result of the previous section, we see that Alexander polynomials are sometimes Salem polynomials (albeit in $-x$).

Indeed, Silver and Williams \cite{SW}, in their study of Mahler measures of Alexander polynomials, found families of links whose Alexander polynomials had Mahler measure equal to a Salem number. The first family $l(q)$ was obtained \cite[Example 5.1]{SW} from the link $7^2_1$ by giving $q$ full right-handed twists to one of the components as it passed through the other component (the trivial knot). The Mahler measure of the Alexander polynomials of these links produced a decreasing sequence of Salem numbers for $q=1,2,\dots,11$. For $q=10$ the Salem number 1.18836\dots (the second-smallest known) was produced, with minimal polynomial$$x^{18}-x^{17}+x^{16}-x^{15}-x^{12}+x^{11}-x^{10}+x^{9}-x^{8}+x^{7}-x^{6}-x^{3}+x^{2}-x+1,$$ while  $q=11$ gave the Salem number $M(L(x))=1.17628\dots$. For $q>11$ Salem numbers were not produced.
The second example was obtained in a similar way \cite[Example 5.8]{SW}, using the link formed from the knot $5_1$ by an adding the trivial knot encircling two strands of the knot, and then giving these strands $q$ full right-hand twists. For increasing $q\ge 3$ this gave a monotonically increasing sequence of Salem numbers tending to the smallest Pisot number $\theta_0=1.3247\dots$. These Salem numbers are equal to the Mahler measure $M(x^{2(q+1)}(x^3-x-1)+x^3+x^2-1)$. Furthermore,  $M(x^{n}(x^3-x-1)+x^3+x^2-1)$ is also a Salem number for $n\ge 9$ and odd. Silver (private communication)  has shown  that these Salem numbers are also Mahler measures of Alexander polynomials: ``Putting an odd number of half-twists in the rightmost arm of the pretzel knot produces 2-component links rather than knots.
Their Alexander polynomials have two variables. However, setting the two variables equal to each other produces the so-called 1-variable
Alexander polynomials, and indeed the `odd' sequence of Salem polynomials \dots  results.''

\subsection{Lengths of closed geodesics}\label{s3} It is known that there is a bijection between the set of Salem numbers and the set of closed geodesics on certain arithmetic hyperbolic surfaces. Specifically, the length of the geodesic is  $2\log\tau$, where $\tau$ is the Salem number corresponding to the geodesic. Thus  there is a smallest Salem number iff there is a geodesic of minimal length among all closed geodesics on all arithmetic hyperbolic surfaces. See Ghate and Hironaka \cite[Section 3.4]{GH} and also Maclachlan and Reid \cite[Section 12.3]{MR} for details.

\subsection{Arithmetic Fuchsian groups} Neumann and Reid \cite[Lemmas 4.9, 4.10]{NR} have shown that  Salem numbers are precisely the spectral radii of hyperbolic elements of arithmetic Fuchsian groups. See also \cite{GH}, \cite[pp. 378--380]{MR} and \cite[Theorem 9.7]{Lei}.

The following result is related.
\bt[{{ Sury \cite{Su} }}] The set of Salem numbers is bounded away from $1$ iff there is some neighbourhood $U$ of the identity in $\SL_2(\R)$ such that, for each arithmetic cocompact Fuchsian group $\Gamma$, the set $\Gamma\cap U$ consists only of elements of finite order.
\et
A Fuchsian group is a subgroup $\Gamma$ of $\PSL_2(\R)$  that acts discontinuously on the upper half-plane $\mathbb H$ (i.e., for $z\in\mathbb H$ no orbit $\Gamma z$ has an accumulation point).

\subsection{Dynamical systems}
\subsubsection{} For given $\beta>1$, define the map $T_\beta:[0,1]\to[0,1)$ by $T_\beta x=\{\beta x\}$, the fractional part of $\beta x$. Then from $x=\frac{\lfloor \beta x\rfloor}{\beta}+\frac{T_\beta x}{\beta}$ we obtain the identity $x=\sum_{n=1}^\infty \frac{\lfloor \beta T_\beta^{n-1} x\rfloor}{\beta^n}$, the (greedy) {\it $\beta$-expansion of $x$} \cite{Pa}.

Klaus Schmidt \cite{Sch} showed that if the orbit of $1$ is eventually periodic for all $x\in\Q\cap[0,1)$ then $\beta$ is a Salem or Pisot number. He also conjectured that, conversely, for $\beta$ a Salem number, the orbit of $1$ is eventually periodic. This conjecture was proved by Boyd \cite{B4} to hold for Salem numbers of degree $4$. However, using a heuristic model in \cite{B5}, his results indicated that while Schmidt's conjecture was likely to also hold for Salem numbers of degree $6$, it may be false for a positive proportion of Salem numbers of degree $8$.  
As Boyd points out, the basic reason seems to be that, for $\beta$ a Salem number of degree $d$, this orbit corresponds to a pseudorandom walk on a $d$-dimensional lattice. Under this model, but assuming true randomness, the probability of the walk intersecting itself is $1$ for $d\le 6$, but is less than $1$ for $d>6$.
Recently, computational degree-8 evidence relating to Boyd's model was compiled by Hichri \cite{Hic}, \cite{Hic2}.

  Hare and Tweddle \cite[Theorem 8]{HT} give examples of Pisot numbers for which the sequences of Salem numbers from Salem's construction  that tend to the Pisot number from above have eventually periodic orbits.
See also \cite{B4.5}. 

For a survey of $\beta$-expansions when $\beta$ is a Pisot or Salem number, see \cite{Ha-conf}. See also Vaz,  Martins Rodrigues, and Sousa Ramos \cite{VMP}.

\subsubsection{} Lindenstrauss and Schmidt \cite [Theorem 6.3]{LS} showed that there  exists ``a connection between two-sided beta-shifts of Salem numbers and the nonhyperbolic ergodic toral automorphisms defined by
 the companion matrices of their minimal polynomials. However, this connection is
much more complicated and tenuous than in the Pisot case'' -- see \cite{Sch2} and \cite[Theorem 6.1]{LS}. 

\subsection{Surface automorphisms }

A {\it K3 surface} is a simply-connected compact complex surface $X$ with trivial canonical bundle. The intersection form  on $H^2(X,\Z)$ makes it into an even unimodular 22-dimensional lattice of signature $(3,19)$; see  \cite[p.17]{mcm5}. Now let $F:X\to X$ be an automorphism of positive entropy of a K3 surface $X$. Then McMullen \cite[Theorem 3.2]{mcm1} has proved that the spectral radius $\la(F)$ (modulus of the largest eigenvalue) of $F$ acting by pullback on this lattice is a Salem number. More specifically,
the characteristic polynomial $\chi(F)$  of the induced map $F^*|H^2(X,\R)\to H^2(X,\R)$ is  the minimal polynomial of a Salem number
 multiplied by $k\ge 0$ cyclotomic polynomials.  Since $\chi(F)$ has degree $22$, the degree of $\la(F)$ is at most 22.  (If $X$ is projective, $X$ has Picard group of rank at most 20,  and so $\la(F)$ has degree at most 20.)

 It is an interesting problem to describe which Salem numbers arise in this way.
McMullen \cite{mcm1} found 10 Salem numbers of degree 22 and trace $-1$, also having some other properties, from which he was able to construct from each of these Salem numbers a K3 surface automorphism having a Siegel disc. (These were the first known examples having Siegel discs.) 
Gross and McMullen \cite{GM} have shown that if the minimal polynomial $S(x)$ of a Salem number of degree 22 has $|S(-1)|=|S(1)|=1$ (which they call the {\it unramified} case) then it is the characteristic polynomial of an automorphism of some (non-projective) K3 surface $X$.
(If the entropy of $F$ is $0$ then this characteristic polynomial is simply a product of cyclotomic polynomials.)
It is known (see \cite [p.211]{mcm1} and references given there) that the topological entropy $h(F)$ of $F$  is equal to $\log\la(F)$, so is either $0$ or the logarithm of a Salem number.

For each even  $d\ge 2$ let $\t_d$ be the smallest Salem number of degree $d$.
McMullen \cite[Theorem 1.2]{mcm3} has proved that if $F : X \to X$ is an automorphism of {\it any} compact complex surface $X$ with
positive entropy, then $h(F) \ge  \log \t_{10}=\log(1.176\dots)=0.162\dots$.
 Bedford and Kim \cite{BK} have shown that this lower bound is realised by a particular rational surface automorphism. McMullen
\cite{mcm4} showed that it was realised for a non-projective K3 surface automorphism, and later \cite{mcm5} that it was realised for a projective K3 surface automorphism.  He showed that the value $\log\t_d$ was realised for a projective K3 surface automorphism
for $d=2, 4, 6, 8, 10$ or $18$, but not for $d=14$, $16$, or $20$. 
(The case $d=12$ is currently undecided.)

Oguiso \cite{O2} remarked that, as for K3 surfaces (see above), the characteristic polynomial of an automorphism of arbitrary compact  K\" ahler surface is also the minimal polynomial of a Salem number  multiplied by $k\ge 0$ cyclotomic polynomials. This is because McMullen's proof for K3 surfaces in \cite{mcm1} is readily generalised. In another paper \cite{O1} he proved that this result also held for automorphisms of hyper-K\" ahler manifolds.  Oguiso \cite{O2} also constructed an automorphism $F$ of a (projective) K3 surface
 with $\la(F)=\t_{14}$.  Here the K3 surface was projective, contained an $E_8$ configuration of rational curves, and the  automorphism also had a Siegel disc.

  Reschke \cite{R} studied the  automorphisms of two-dimensional complex tori. He showed that the entropy  of such an automorphism, if positive,  must be a Salem number of degree at most 6, and gave necessary and sufficient conditions for such a Salem number to arise in this way.

\subsection{Salem numbers and Coxeter systems} Consider a Coxeter system $(W,S)$, consisting of a multiplicative group $W$ generated by a finite set $S=\{s_1,\dots,s_n\}$, with relations $(s_is_j)^{m_{ij}}=1$ for each $i,j,$ where $m_{ii}=1$ and $m_{ij}\ge 2$ for $i\ne j$.  The $s_i$ act as reflections on $\mathbb R^n$. 
 For any $w\in W$ let $\la(w)$ denote its spectral radius. This is the modulus of the largest eigenvalue of its action on $\R^n$.
Then  McMullen \cite[Theorem 1.1]{mcm2} shows that when $\la(w)>1$ then $\la(w)\ge \t_{10}=1.176\dots$. This could be interpreted as circumstantial evidence for $\t_{10}$ indeed being the smallest Salem number.

The Coxeter diagram  of $(W, S)$ is the weighted graph whose vertices are the
set $S$, and whose edges of weight $m_{ij}$ join $s_i$ to $s_j$ when $m_{ij} \ge 3$.
Denoting by $Y_{a,b,c}$  the Coxeter system whose diagram is a tree with 3 branches of
lengths $a$, $b$ and $c$, joined at a single node, McMullen also showed that
the smallest Salem numbers of degrees 6, 8 and 10 coincide with $\la(w)$ for
the Coxeter elements of $Y_{3,3,4}$, $Y_{2,4,5}$ and $Y_{2,3,7}$ respectively.
 In particular,  $\la(w)=\t_{10}$ for the Coxeter elements of $Y_{2,3,7}$.

\subsection{Dilatation of pseudo-Anosov automorphisms} For a closed connected oriented surface $\mathcal S$ having a pseudo-Anosov automorphism  that is a product of two positive multi-twists, Leininger \cite[Theorem 6.2]{Lei}  showed that its dilatation is at least $\t_{10}$.  This follows from McMullen's work on Coxeter systems quoted above. The case of equality is explicitly described (in particular, $\mathcal S$ has genus 5). 
(However, on  surfaces of genus $g$  there are examples of pseudo-Anosov automorphisms  having  dilatations equal to $1+O(1/g)$ as $g\to\infty$. These are not Salem numbers when $g$ is sufficiently large.)

\subsection{Bernoulli convolutions} Following Solomyak \cite{So}, let $\la\in(0,1)$, and $Y_\la=\sum_{n=0}^\infty\pm \la^n$, with the $\pm$ chosen independently `$+$' or `$-$' each with probability $\frac12$. Let $\nu_\la(E)$ be the probability that $Y_\la\in E$, for any Borel set $E$. So it is the infinite convolution product of the means $\frac12(\delta_{-\la^n}+\delta_{\la^n})$ for $n=0,1,2,\dots,\infty$, and so is called a {\it Bernoulli convolution}. Then $\nu_\la(E)$ satisfies the self-similarity property
\[
\nu_\la(E)=\frac12\left(\nu_\la(S_1^{-1}E)+\nu_\la(S_2^{-1}E)\right),
\]
where $S_1x=1+\la x$ and $S_2x=1-\la x$. It is known that the support of $\nu_\la$ is a Cantor set of zero length when $\la\in(0,\frac12)$, and the interval $[-(1-\la)^{-1},(1-\la)^{-1}]$ when $\la\in(\frac12,1)$. When $\la=\frac12$, $\nu_\la$ is the uniform measure on $[-2,2]$. Now the Fourier transform $\hat\nu_\la(\xi)$ of $\nu_\la$ is equal to $\prod_{n=0}^\infty \cos(\la^n\xi)$. Salem \cite[p. 40]{Sam} proved that if $\la\in(0,1)$ and $1/\la$ is not a Pisot number, then $\lim_{\xi\to\infty}\hat\nu_\la(\xi)=0$. This contrasts with an earlier result of Erd\H os that if $\la\ne \frac12$ and $1/\la$ is a Pisot number, then $\hat\nu_\la(\xi)$ does not tend to $0$ as $\xi\to\infty$. Furthermore Kahane \cite{Ka} (as corrected in  \cite[p. 10]{PSS}) proved that if $1/\la$ is a Salem number then for each $\eps>0$
\begin{equation}\label{E-xi}
\limsup_{\xi\to\infty} |\hat\nu_\la(\xi)|\,|\xi|^\eps=\infty.
\end{equation}
The only property of Salem numbers used in the proof of \eqref{E-xi} is Proposition \ref{P-Teps}.
Now consider the set $S_{a,\gamma}$ of all $\la\in(a,1)$ such that $\int_{-\infty}^\infty |\hat\nu_\la(\xi)|^2\,|\xi|^{2\gamma} d\xi<\infty$. Then Peres,  Schlag and  Solomyak \cite[Proposition 5.1]{PSS}
use this set, Lemma \ref{L-2} and \eqref{E-xi}, to give a sufficient condition for the set $T$ of Salem numbers to be bounded away from $1$. Specifically, they prove that if there exist some $\gamma>0$ and $a<1$ such that $S_{a,\gamma}$ is a so-called {\it residual set} (i.e., it is the intersection of countably many sets with dense interiors), then $1$ is not a limit point of $T$.

Concerning the behaviour of $\hat\nu_\la(\xi)$ as $\xi\to\infty$, where $1/\la=\theta>1$ is an algebraic integer: if $\theta$ has another conjugate $\theta'\ne \theta$ outside the unit circle, then Bufetov and  Solomyak \cite[Corollary A.3]{BS} recently proved that $\hat\nu_\la(\xi)$ is bounded by a negative power of $\log\xi$ as $\xi\to\infty$. 
On the other hand, it is a result of Erd\H os \cite{Er} that when $\theta$ is a Pisot number then $\limsup_{\xi\to\infty}\hat\nu_\la(\xi)$  is positive, showing that such a bound is not possible for Pisot numbers. Bufetov and  Solomyak ask, however, whether a bound of that kind might hold for Salem numbers (the only remaining undecided case for an algebraic integer $\theta>1$).

Recently Feng \cite{Fe} has also studied $\nu_\la$ when $1/\la$ is a Salem number, proving that then the corresponding measure $\nu_\la$ is a multifractal measure satisfying the multifractal formalism in all of the increasing part of its multifractal spectrum.

For two very readable surveys of Bernoulli convolutions, including connections with Salem numbers, see Kahane's `Reflections on Paul Erd\H os \dots' article \cite{K-AMS} and the much longer survey by Peres,  Schlag and  Solomyak \cite{PSS} referred to above.

{\bf Acknowledgements.} This paper is an expanded version of a talk
that I gave at  the meeting `Growth and Mahler measures in geometry
and topology' at the Mittag-Leffler Institute, Djursholm, Sweden,
in July 2013. I would like to thank the organisers Eriko Hironaka
and Ruth Kellerhals for the invitation to attend the meeting, and to
thank them, the Institute staff and fellow participants for making it such a stimulating
and pleasant week.

I also thank David Boyd, Yann Bugeaud,  Nigel Byott,  Eriko Hironaka, Aleksander Kolpakov, James McKee, Curtis McMullen, Andrew Ranicki, Georges Rhin, Dan Silver, Joe Silverman, Boris Solomyak, Timothy Trudgian   and especially the referee for their very helpful comments and corrections on  earlier drafts of this survey.

\end{document}